\documentclass{mathincs}
\usepackage[T1]{fontenc}
\usepackage{floatflt}
\usepackage{algorithm,algpseudocode}
\usepackage{graphicx,caption}
\usepackage{amsfonts}
\usepackage[unicode, pdftex]{hyperref}
\usepackage{array}
\usepackage{longtable}
\usepackage{multirow}
\usepackage{vcell}

 \newtheorem{thm}{Theorem}[section]
 \newtheorem{theorem}[thm]{Theorem}
 \newtheorem{example}[thm]{Example}

 \newtheorem{conjecture}[thm]{Conjecture}
 \theoremstyle{definition}
 \newtheorem{definition}[thm]{Definition}
 \theoremstyle{remark}

 \numberwithin{equation}{section}

\begin{document}
\title{A review on computational aspects of polynomial amoebas}
\author{Vitaly A. Krasikov}

\address{%
Plekhanov Russian University of Economics \\
36 Stremyanny Lane\\
Moscow, Russia}
%\orcidID{0000-0001-9686-8493}
%\authorrunning{Vitaly A. Krasikov}
%\institute{Plekhanov Russian University of Economics, Moscow, Russia \\
\email{krasikov.va@rea.ru,vitkras@inbox.ru}
\thanks{The research was supported by a grant from the Russian Science Foundation No.22-21-00556, \url{https://rscf.ru/project/22-21-00556/}}


\maketitle              % typeset the header of the contribution
\begin{abstract}
We review results of papers written on the topic of polynomial amoebas with an emphasis on computational aspects of the topic. The polynomial amoebas have a lot of applications in various domains of science. Computation of the amoeba for a given polynomial and describing its properties is in general a problem of formidable complexity. We describe the main algorithms for computing and depicting the amoebas and geometrical objects associated with them, such as contours and spines. We review the latest software packages for computing the polynomial amoebas and compare their functionality and performance.

%\keywords{polynomial amoebas \and Newton polytope \and visualization algorithms.}
\end{abstract}

\section{Introduction}

An amoeba of a polynomial is the projection of its zero locus on the space of absolute values in a logarithmic scale. Initially this notion was introduced by I.M. Gelfand, M.M. Kapranov, and A.V.~Zelevinsky in 1994~\cite{GKZ:1994}. The polynomial amoebas have numerous applications in topology~\cite{Guilloux-Marche:2021}, dynamical systems~\cite{Theobald:2002}, algebraic geometry~\cite{Mikhalkin:2017,Lang:2019,Goucha-Gouveia:2021}, complex analysis~\cite{Forsberg-Passare-Tsikh:2000}, mirror symmetry~\cite{Ruan:2000,Feng-He-Kennaway-Vafa:2008}, measure theory~\cite{Mikhalkin-Rullgard:2001,Passare-Rullgard:2000}, statistical physics~\cite{Passare-Pochekutov-Tsikh:2013}.

There are a lot of surveys published on the topic of polynomial amoebas (see~\cite{Mikhalkin:2004,Yger:2012,deWolff:2017}), but the present paper is focused on its computational aspects. The computation of amoebas is a~problem of great importance, but it also has a formidable complexity for the following reasons. First, an amoeba of a polynomial is in general an unbounded set in~$\mathbb{R}^n$ but at the same time in higher dimensions the amoeba for any given polynomial is very <<narrow>> subset of~$\mathbb{R}^n,$ so even locating the amoeba in the whole space is the non-trivial problem. Another complex problem is locating connected components of the amoeba complement, since their number can differ and their size can be very small. A possible significant difference between the sizes of the amoeba and connected components of its complement is the reason why the computation for complex amoebas requires very high precision.~\cite{Zhukov-Sadykov:2023}

The general complexity of the amoeba computation relates to some partial problems as well. For example, the problem of deciding whether a given point belongs to the amoeba (or the membership problem) in general cannot be solved in polynomial time. Since the membership problem is crucial for computing the amoebas, there are a lot of papers containing algorithms for a fast solution to it~\cite{Theobald:2002,Forsgard-Matusevich-Mehlhop-deWolff:2018,Purbhoo:2008,Theobald-deWolff:2015}).

The main purpose of the present survey is to take a wide look at the concept of polynomial amoeba in the computational context, including some historical details on the origin of polynomial amoebas, their modern applications and computational aspects. It includes the review on software packages for computing the amoebas and benchmarking results for these packages. 

\section{Basic notation and motivation}

Let $p (z)$ be a (Laurent) polynomial in $n$ complex variables:
$$p(z_1,\ldots, z_n)=\sum\limits_{\alpha\in A} c_{\alpha} z^{\alpha}=\sum\limits_{\alpha\in A} c_{\alpha_1\ldots\alpha_n}z_1^{\alpha_1}\cdot\ldots\cdot z_n^{\alpha_n},$$
where $A \subset \mathbb{Z}^n$ is a finite set.

\begin{definition} \rm
{\it The amoeba} $\mathcal{A}_{p}$ of a polynomial $p(z)$ is the image of its zero locus under the map $\textrm{Log}: \left(\mathbb{C^*}\right)^n \rightarrow \mathbb{R}^n,$ where $\mathbb{C^*}=\mathbb{C}\backslash \{0\}:$ $$\textrm{Log}: (z_1, \ldots, z_n) \longmapsto (\textrm{ln}|z_1|, \ldots, \textrm{ln}|z_n|).$$
\end{definition}

Initially this notion was introduced by I.M. Gelfand, M.M. Kapranov, and A.V. Zelevinsky in~1994~\cite{GKZ:1994}. The resemblance between the shapes of the image Log $({z|p(z) = 0})$ for $n = 2$ and the unicellular organism is the reason for using the biological term.

In the similar way the notion of amoeba can be introduced for algebraic or even transcendental hypersurfaces~\cite{Passare-Pochekutov-Tsikh:2013} but in what follows the focus is on the polynomial amoebas.

In addition to the Log map one may consider also the map $\textrm{Arg}:  \left(\mathbb{C^*}\right)^n \rightarrow \left(S^1\right)^n:$

$$\textrm{Arg}: (z_1, \ldots, z_n) \longmapsto (\textrm{arg}(z_1), \ldots, \textrm{arg}(z_n)).$$ 

\begin{definition} \rm
{\it The coamoeba} $co\mathcal{A}_{p}$ of a polynomial $p(z)$ is the image of its zero locus under the map $\textrm{Arg}.$
\end{definition}
In some sources (see, for example,~\cite{Feng-He-Kennaway-Vafa:2008}) this image is also called <<alga>> of the hypersurface. Coamoebas were introduced by M. Passare in a talk in 2004. 

There is a strong connection between characteristics of the amoeba for a polynomial and its Newton polytope.

\begin{definition} \rm
The convex hull in~$\mathbb{R}^n$ of the set~$A$ is called {\it the Newton polytope of} $p(z),$ we denote it by~$\mathcal{N}_p.$ 
\end{definition}

The connected components of the amoeba complement $\textrm{}^c\mathcal{A}_p = \mathbb{R}^n\backslash \mathcal{A}_p$ are convex subsets in~$\mathbb{R}^n.$ They are in bijective correspondence with the different Laurent expansions (centered at the origin) of the rational function~$1/p(z)$~\cite{GKZ:1994}.

\begin{definition} \rm
{\it The compactified amoeba} $\bar{\mathcal{A}}_{p}$ of a polynomial $p(z)$ is the closure of the image of its zero locus under the map $\nu: (z_1, \ldots, z_n) \longmapsto \frac{\sum_{\alpha\in A}|z^\alpha|\cdot\alpha}{\sum_{\alpha\in A}|z^\alpha|}.$
\end{definition}

In general case, the following statement holds:

\begin{theorem} \rm~\cite{Forsberg-Passare-Tsikh:2000}
The number of connected components of the amoeba complement $\textrm{}^c\mathcal{A}_p$ is at least equal to the number of vertices of the Newton polytope $\mathcal{N}_p$ and at most equal to the total number of integer points in $\mathcal{N}_p\cap \mathbb{Z}^n.$
\label{thm:connectedComponentsNumber}
\end{theorem}

The part of Theorem~\ref{thm:connectedComponentsNumber} concerning lower bound for the number of connected components was proved in~\cite{GKZ:1994} and~\cite{Mkrtchian-Yuzhakov:1982}.

\begin{definition} \rm
If the number of connected components of $\textrm{}^c\mathcal{A}_p$ is equal to the number of vertices of the Newton polytope~$\mathcal{N}_p,$ the amoeba~$\mathcal{A}_p$ is called {\it solid.} If the number of connected components of $\textrm{}^c\mathcal{A}_p$ is equal to the number of integer points in~$\mathcal{N}_p\cap \mathbb{Z}^n,$ the amoeba~$\mathcal{A}_p$ is called {\it optimal.}
\end{definition}

\begin{example} \rm
Consider a family of polynomials $V_\alpha$ defined through the generating function (see~\cite[Chapter 2.3]{Xu-Dunkl:2014}) $$(1-2\langle a,x\rangle+||a||^2)^\frac{n-1}{2}=\sum\limits_{\alpha\in\mathbb{N}^n_0}a^\alpha V_\alpha(x).$$

By neglecting a monomial factor of $V_\alpha$ and replacing variables $x_i, i=1,\ldots,n$ by $\xi_i=x_i^2,$ one may obtain polynomial~$\tilde V_\alpha(\xi)$ with the same number of components in its amoeba complement as $V(x)$ (see Lemma 2.6 and Example 2.10 in~\cite{Bogdanov-Sadykov:2020}).

The Newton polytope and the amoeba of the polynomial $\tilde V_{(5,2,3)}(\xi)$ are shown in Figure~\ref{fig:amoeba_3d}. Components of the amoeba complement are depicted there as colored convex shapes, the amoeba itself is a white space in between. The Newton polytope for~$\tilde V_{(5,2,3)}(\xi)$ contains~$12$ integer points and this number coincides with the number of connected components in the amoeba complement, so this amoeba is optimal.

\begin{figure}[h!]
\hskip0.5cm
\begin{minipage}{5.5cm}
\includegraphics[width=5.5cm]{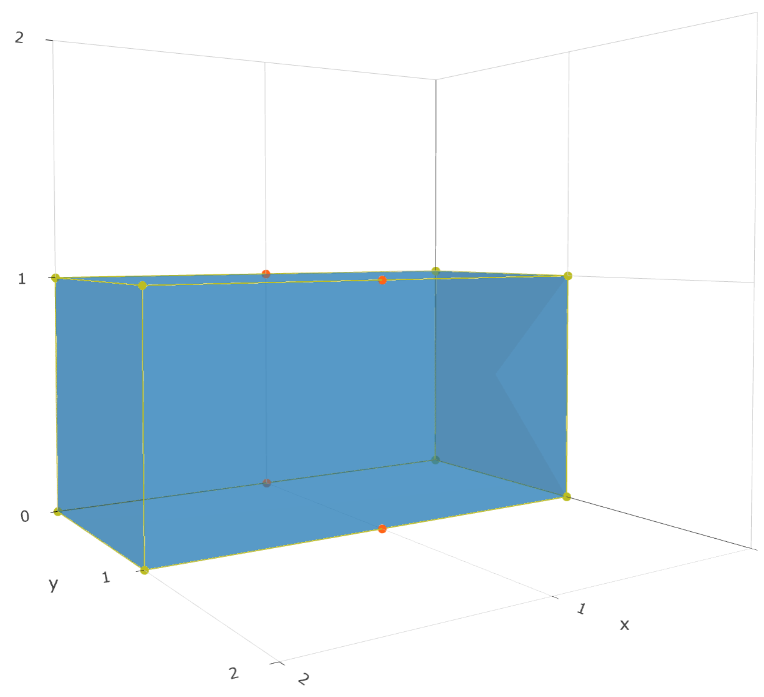}
%\centering $(a)$
\end{minipage}
%\hskip1cm
\begin{minipage}{7.5cm}
\includegraphics[width=7.5cm]{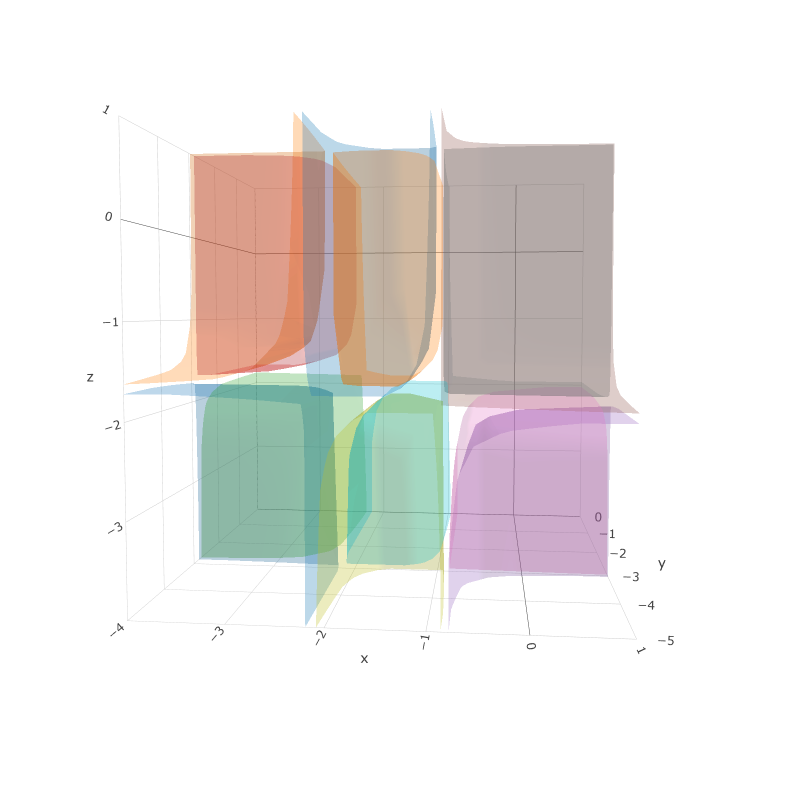}
%\centering $(b)$
\end{minipage}

\caption{The Newton polytope and the amoeba of the polynomial $\tilde V_{(5,2,3)}$}
\label{fig:amoeba_3d}
\end{figure}

\end{example}

It follows from Theorem~\ref{thm:connectedComponentsNumber} that there exists an injective function from the set of connected components of $\textrm{}^c\mathcal{A}_p$ to $\mathcal{N}_p\cap \mathbb{Z}^n.$ One may construct such a function by using the Ronkin function $$N_p(x)=\left(\frac{1}{2\pi i}\right)^n \int_{\textrm{Log}^{-1}(x)}\textrm{ln}|p(z)|\frac{dz}{z},\quad x\in\mathbb{R}^n.$$

The function $N_p$ is convex on $\mathbb{R}^n$ and it is affine linear on an open connected set $\Omega\subset\mathbb{R}^n$ if and only if $\Omega\subset\textrm{}^c\mathcal{A}_p.$ 

An analogue of the Ronkin function for coamoebas is introduced in~\cite{Johansson-Samuelsson:2017}.

\begin{definition} \rm~\cite{Forsberg-Passare-Tsikh:2000}
{\it The order of the connected component} $E\subset\textrm{}^c\mathcal{A}_p$ is defined as vector $\nu\in\mathbb{Z}^n$ with the components $$\nu_j=\frac{1}{(2\pi i)^n}\int_{\textrm{\scriptsize Log}^{-1}(u)}\frac{z_j\partial_j f(z)}{f(z)}\frac{dz_1\wedge\ldots\wedge dz_n}{z_1\ldots z_n},\quad j=1,\ldots,n,$$ where $u\in E.$ This number does not depend on the choice of~$u,$ since the homology class of the cycle~$\textrm{Log}^{-1}(u)$ is the same for all $u\in E.$
\end{definition}

The order of the component induces injective function $$\nu: \{E\} \rightarrow \mathbb{Z}_n\cap\mathcal{N}_p.$$

%\begin{definition} \rm
%{\it The recession cone} of a convex set $M$ is the set-theoretical maximal element in the family of convex cones whose shifts are contained in~$M.$
%\end{definition}
%
%\begin{definition} \rm
%For any $\nu\in\mathbb{R}^n\cap\mathcal{N}_p$ {\it the dual cone} $C_\nu$ is a set $$C_\nu=\{s\in\mathbb{R}^n| <s,\nu> = \max\limits_{\alpha\in\mathcal{N}_p} <s,\alpha> \}.$$
%\end{definition}
%
%If the component~$E$ of the complement~$\textrm{}^c\mathcal{A}_p$ has order~$\nu\in\mathcal{N}_p,$ then the dual cone~$C_\nu$ equals to the recession cone of~$E.$

%There are some geometric notions associated with an amoeba.

\begin{definition} \rm
{\it The contour} of the amoeba~$\mathcal{A}_{p}$ is the set~$\mathcal{C}_{p}$ of critical points of the logarithmic map Log restricted to the zero locus of the polynomial $p(x).$ 
\end{definition}

The contour is the closed real-analytic hypersurface in~$\mathbb{R}^n, $ the boundary $\partial\mathcal{A}_p$ is a subset of the contour~$\mathcal{C}_p$ but is in general different from it. Results on the maximal number of intersection points of a line with the contour of hypersurface amoebas are given in~\cite{Lang-Shapiro-Shustin:2021}.

Let $\mathcal{A}'\subset\mathbb{R}^n\cap\mathcal{N}_p$ be a set of vectors $\alpha$ such that $\textrm{}^c\mathcal{A}_p$ contains components of the order~$\alpha$ and $$a_\alpha=\frac{1}{(2\pi i)^n}\int_{\textrm{\scriptsize Log}^{-1}(u)}\textrm{Log}|\frac{f(z)}{z^\alpha}|\frac{dz_1\wedge\ldots\wedge dz_n}{z_1\ldots z_n},\quad\forall\alpha\in\mathcal{A}',\quad u\in E_\alpha.$$

Consider the function $S_p:\mathbb{R}^n\rightarrow\mathbb{R}:$

$$S_p(x)=\max\limits_{\alpha\in\mathcal{A}'}\{<\alpha,x>+a_\alpha\}$$

\begin{definition} \rm
The set $\left\{x\in\mathbb{R}^n,\textrm{ where }S_p(x)\textrm{ is nonsmooth}\right\}$ is called {\it the spine of}~$\mathcal{A}_p.$
\end{definition}

\begin{example} \rm
Consider the polynomial $p_1(z_1,z_2)=5 z_1 + 15 z_1^2 + 8 z_1^3 z_2 + 10 z_2^2 + 10 z_1^3 z_2^2 + 8 z_2^3 + 15 z_1 z_2^4 + 5 z_1^2 z_2^4 + 50 z_1 z_2^3 + 50 z_1^2 z_2.$ The Newton polytope, the amoeba, its spine, and the compactified amoeba of~$p_1$ are shown in Figure~\ref{fig:p1_NewtonPolytope}.

\begin{figure}[h!]
\hskip0.5cm
\begin{minipage}{4cm}
\includegraphics[width=3.5cm]{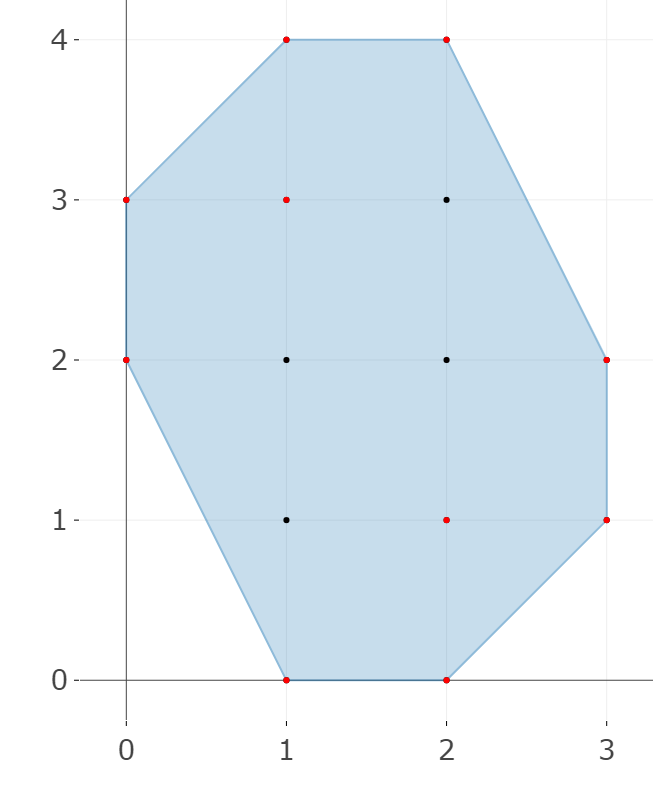}
\centering $(a)$
\end{minipage}
\hskip1cm
\begin{minipage}{5cm}
\includegraphics[width=5cm]{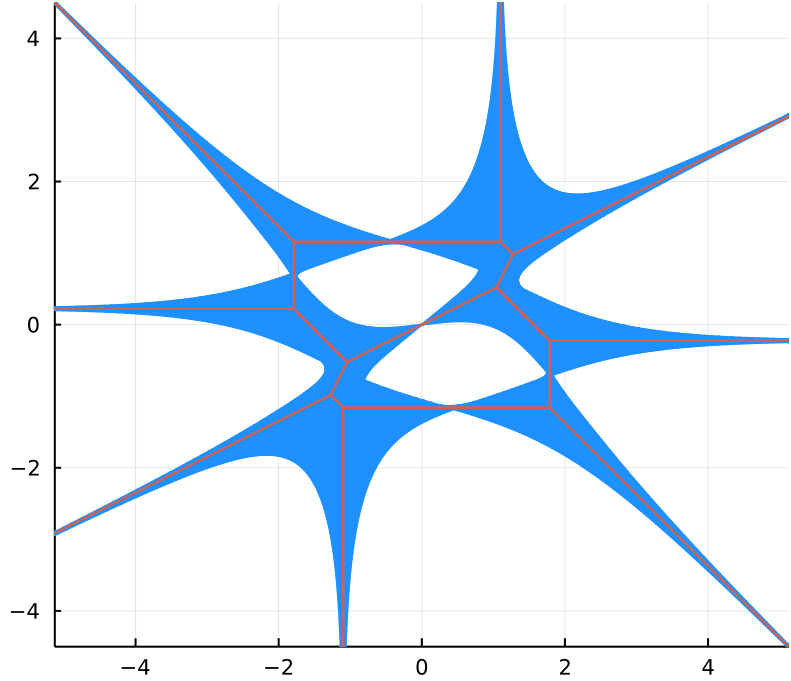}
\centering $(b)$
\end{minipage}

\begin{minipage}{5.5cm}
\includegraphics[width=5.5cm]{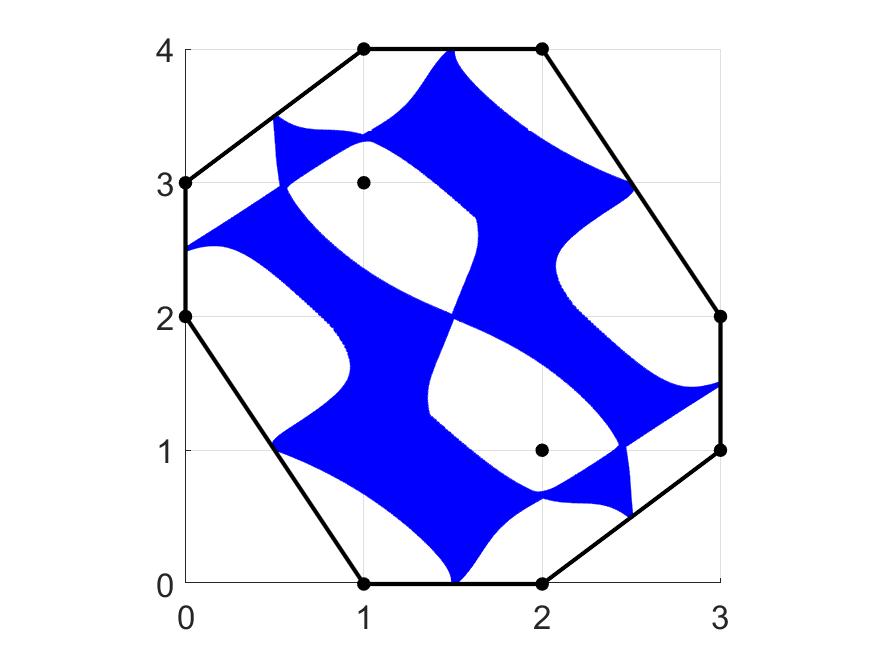}
\centering $(c)$
\end{minipage}
\begin{minipage}{5cm}
\includegraphics[width=5cm]{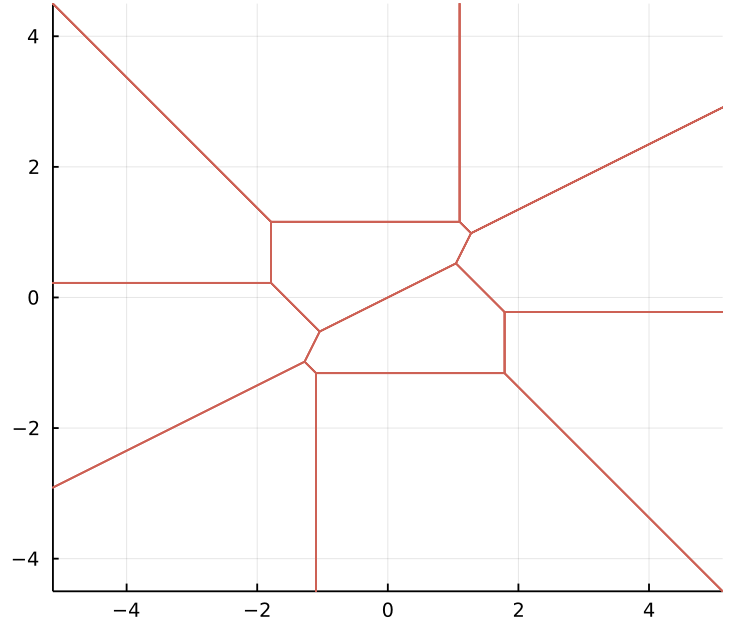}
\centering $(d)$
\end{minipage}
\caption{The Newton polytope (a), the amoeba (b), the compactified amoeba (c) and the spine of the amoeba (d) of the polynomial $p_1.$}
\label{fig:p1_NewtonPolytope}
\end{figure}
\end{example}

Motivation for this survey is the research project currently conducted by scientific group in the Laboratory of artificial intelligence, neurotechnology and business analytics in Russian Plekhanov University of Economics. Main problem under consideration is the following conjecture by M. Passare~\cite{PassareConjectureFormulation:2009}: 

\begin{conjecture} (Passare conjecture) \rm
Let $p(z)$ be a maximally sparse polynomial (i.e. the support of $p$ is equal to the set of vertices of its Newton polytope~$\mathcal{N}_p$).  Then the amoeba~$\mathcal{A}_p$ is solid.
\end{conjecture}

M. Passare and H. Rullgård proved that if the number of vertices is less than or equal to $n + 2$, then the spine is contained in the amoeba~\cite{Passare-Rullgard:2000}. M. Passare stated that <<it would seem very plausible that the number of complement components is minimal for maximally sparse polynomials with at most $n + 2$ terms>>.

\begin{example} \rm
Consider polynomial~$\tilde{p}_1(z_1, z_2),$ obtained by dropping from~$p_1(z_1, z_2)$ monomials corresponding to inner lattice points of its Newton polytope: $\tilde{p}_1(z_1, z_2) = 5 z_1 + 15 z_1^2 + 8 z_1^3 z_2 + 10 z_2^2 + 10 z_1^3 z_2^2 + 8 z_2^3 + 15 z_1 z_2^4 + 5 z_1^2 z_2^4.$ Amoeba of $\tilde{p}_1$ is shown in Figure~\ref{fig:p1_solid_Amoeba}.

\begin{figure}[h!]
%\hskip0.5cm
\centering
\begin{minipage}{5cm}
\includegraphics[width=5cm]{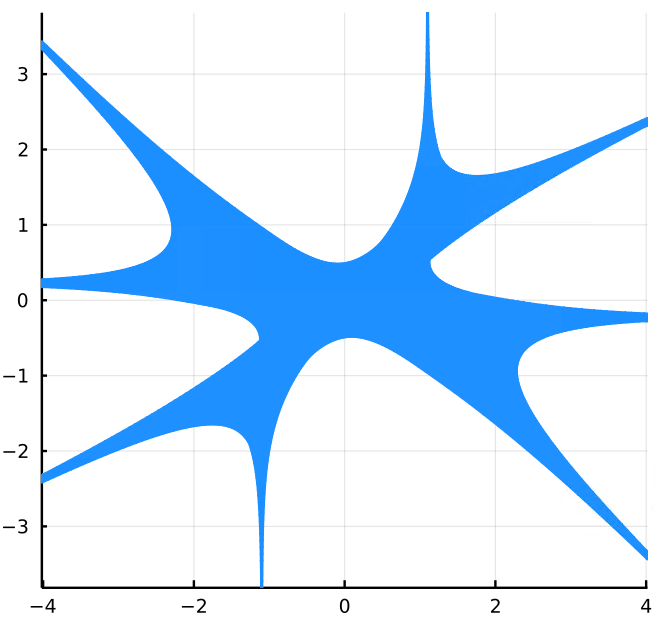}
\end{minipage}

\caption{The amoeba of the polynomial $\tilde{p}_1.$}
\label{fig:p1_solid_Amoeba}
\end{figure}
\end{example}

There are results on solid amoebas for some classes of sparse polynomials. In~\cite{Iliman-deWolff:2016} authors consider the class of polynomials of several real variables, whose Newton polytopes are simplices and the supports are given by all the vertices of the simplices and one additional interior lattice point in the simplices. It is proved that under some conditions amoebas for such polynomials are solid.

In the context of the Passare conjecture the theory of fewnomials developed by A.G. Khovanskii~\cite{Khovanskii:1991} should be mentioned. The whole ideology behind it is that a variety defined by systems of <<simple>> polynomials should have a <<simple>> topology. Maximally sparse polynomial is <<simple>> in the same way, that is, it contains minimal possible number of monomials among all of the polynomials with a given Newton polytope.

\section{Historical reference and applications of polynomial amoebas}

Though the notion of polynomial amoeba was introduced in 1994, some of the results connected to it were known long before this. This sections contains reference to some classical problems, which nowadays are connected to the polynomial amoebas as well as the latter works on the topic.

\subsection{Hilbert's 16th problem}

One of the problems related to the amoebas is the famous Hilbert’s 16th problem on the topology of algebraic curves and surfaces~\cite{Hilbert:1902}. Originally
it was formulated as consisting of two parts: first considering the relative positions of the branches of real algebraic curves of degree $n$ and the second one on the upper bound for the number of limit cycles in two-dimensional polynomial vector fields of degree $n$ and their relative positions.

The first part of Hilbert’s 16th problem followed A. Harnack’s invesigation on the number of separate connected components for algebraic curves of degree $n$ in $\mathbb{R}\mathbb{P}^2$. Harnack proved~\cite{Harnack:1876} that the number of components of such curve does not exceed $\frac{(n-1)(n-2)}{2}+1$. A curve for which this maximum is attained and components have the best possible topological configuration is called a {\it Harnack curve}~\cite{Mikhalkin:2000}. The amoeba of a Harnack curve can be described analytically as $$\left\{(x,y)\in\mathbb{R}^2: \prod\limits_{\alpha,\beta=\pm 1} P_W(\alpha e^x,\beta e^y)\leq 0\right\},$$
where $P_W$ is the polynomial defining the curve~\cite{Feng-He-Kennaway-Vafa:2008}.

Amoebas of Harnack curves enjoy different optimal properties: for example, the amoeba of Harnack curve with genus $g$ has exactly $g$ compact components in its complement~\cite{Kenyon-Okounkov-Sheffield:2006}, a curve is Harnack if and only if its amoeba has the maximal possible area for a given Newton polygon~\cite{Mikhalkin-Rullgard:2001}. It is proved in~\cite{Kenyon-Okounkov:2006} that areas of the amoeba holes (that is, bounded components of the amoeba complement) can be chosen as local coordinates on the manifold of Harnack curves with given boundary and genus.

An explicit integral formula is presented in~\cite{Passare:2016} for the amoeba-to-coamoeba mapping in the case of polynomials defining Harnack curves as well as the exact description of the coamoebas of such polynomials.

\subsection{Amoebas and the convergence of power series}

Consider a power series centered at $0:$
\begin{equation}
\sum\limits_{m\in\mathbb{N}^n}c_mz^m, c_m\in\mathbb{C}, z^m=z_1^{m_1}z_2^{m_2}\ldots z_n^{m_n}.
\label{eq:powerSeries}
\end{equation}

The domain of convergence of the series~\ref{eq:powerSeries} is a complete Reinhardt domain centered at $0$. For $n = 1$ it is stated in Abel’s lemma. O. Cauchy formulated this fact in 1831 \cite[vol. 8, p. 151]{Cauchy-OeuvresCompletes} as follows:

{\it $x$ d\'esignant une variable r\'eelle ou imaginaire, une fonction r\'eelle ou imaginaire de $x$ sera d\'eveloppable en une s\'erie convergente ordonn\'ee suivant les puissances ascendantes de celte variable, tant que le module de la variable conservera une valeur inf\'erieure \`a la plus petite de celles pour lesquelles la fonction ou sa d\'eriv\'ee cesse d'\^etre finie ou continue,} or, in translation,

{\it $x$ designating a real or imaginary variable, a real or imaginary function of $x$ can be developed into an ordered convergent series according to the ascending powers of this variable, as long as the modulus of the variable retains a value lower than the smallest of those for which the function or its derivative
ceases to be finite or continuous.}

The proof of this statement in the case of multiple variables is given in a lot of different sources, but it is not an easy task to detect the initial one.

If all of the singularities of the series~\ref{eq:powerSeries} belong to an algebraic hypersurface then the structure of its domain of convergence is closely related to the amoeba of this hypersurface. More precisely, for any polynomial $p(z)$ there exists a bijective correspondence between connected components of the
amoeba complement for an algebraic set $\{p(z) = 0\}$ and domains of convergence of the Laurent series expansions with denominator $p(z)$~\cite{GKZ:1994}.

%{\Large \bf Add an example. }

%\begin{example}~\cite{Forsberg-PHD}
%
%Consider the polynomial~$p(z_1,z_2)=1+z_1+z_2$ in $\mathbb{C}^2_*.$ The complement of its amoeba $\textrm{}^c\mathcal{A}_p = \mathbb{R}^2\backslash \mathcal{A}_p$ consists of three disjoint components. To these three sets there correspond three different Laurent series whose supports of summation are contained in the sectors indicated in Figure {\bf \large add figure}. The components of $\textrm{}^c\mathcal{A}_p$ are the natural logarithmic domains of convergence for the Laurent series and these components, in accordance with Abel's lemma, contains a translate of exactly one of these cones. Figure {\bf \large add figure} shows the amoeba, its complement and the translated dual cones.
%
%\end{example}

\subsection{Modern amoebas}

There are a lot of papers on different topics connected to the polynomial amoebas. Amoebas can be defined for different varieties including spherical tropical varieties\cite{Kaveh-Manon:2019}, non-hypersurface ones~\cite{Juhnke-Kubitzke-deWolff:2016} and subvarieties of non-commutative Lie groups~\cite{Mikhalkin-Shkolnikov:2022}. Zero-dimensional case is considered in~\cite{Nisse:2016}. Results on the dimension of amoebas of varieties are presented in~\cite{Mikhalkin:2017,Draisma-Rau-Yuen:2020}

The special kind of polyhedral complex which is the subset of Newton polytope of a polynomial and enjoys the key topological and combinatorial properties of the polynomial amoeba is considered in~\cite{Nisse-Sadykov:2019}.

In~\cite{Theobald-deWolff:2016} a new approach to the classical problem of defining the behaviour of univariate trinomial roots is proposed, based on the terms of the polynomial amoeba theory.

In~\cite{Lyapin:2009} amoebas are used for investigating the properties of Riordan's arrays which arise as solutions of Cauchy problem for difference equations.

The polynomial amoebas enjoy some optimal properties if they are defined by hypergeometric functions. It is proved in~\cite{Bogdanov-Sadykov:2020} that under certain nondeneracy conditions amoebas of hypergeometric polynomials are optimal. In~\cite{Passare-Sadykov-Tsikh:2005} it is shown that the singular hypersurface of any nonconfluent hypergeometric function has a solid amoeba. For hypergeometric series the description of the convergence domain in terms of the amoeba complement is given in~\cite{Passare-Tsikh:2002}, the latter results on the topic are~\cite{Nilsson-Passare-Tsikh:2019,Cherepanskiy-Tsikh:2020}.

A lot of natural connections to polynomial amoebas lie in the field of tropical geometry. There are a plenty of surveys and articles on the topic, see, for example, \cite{Mikhalkin:2006,Einsiedler-Kapranov-Lind:2006,Yger:2012,Jensen-Leykin-Yu:2016,Jonsson:2016,deWolff:2017,Forsgard:2020,Lang:2020,Hicks:2020,Kim-Nisse:2021}

%Some of the geometric properties of the polynomial amoebas can be better described in terms of tropical geometry.
%
%Basic object for this kind of geometry is the {\it tropical semiring} $(\mathbb{R}\cup \{\infty\}, \oplus, \odot),$ which is basically a set of real numbers with infinity and operations of addition and multiplication being defined as follows:
%
%$$x \oplus y = \textrm{max} (x,y), \quad x\odot y = x + y.$$
%
%Various frames for developing calculus over this structure are described in~\cite{Yger:2012}
%
%In  $(\mathbb{R}\cup \{\infty\}, \oplus, \odot)$ for the variables $x = (x_1,x_2,\ldots,x_n), x_i\in (\mathbb{R}\cup \{\infty\}, \oplus, \odot), i=1,\ldots,n$ we define {\it tropical polynomial}:
%
%$$p(x_1,\ldots,x_n) = \bigoplus\limits_{k=1}^p a_{i^k} \odot x^{\odot i^k} = \textrm{max}\{a_{i_1}+<x, i_1>, \ldots a_{i_p}+<x, i_p>\}.$$
%
%Tropical polynomial in~$n$ variables is a piecewise-linear concave function on $\mathbb{R}^n$ with integer coefficients.
%
%\begin{definition} \rm
%{\it Tropical hypersurface} $\mathcal{H}(p)$ is a set of points~$x\in \mathbb{R}^n,$ where $p(x)$ is not linear function.  
%\end{definition}

\subsection{Applications in physics}

There are two-dimensional combinatorial systems defined in the domain of mathematical physics, called dimer models~\cite{Kenyon:2008}. In some cases, dimer model graphs for curves in $\mathbb{C}^3$ agree with amoebas for these curves~\cite{Feng-He-Kennaway-Vafa:2008}. In~\cite{Kenyon-Okounkov-Sheffield:2006} the spectral curve of the Kasteleyn operator of the graph is considered. The amoeba of the spectral curve represents the phase diagram of the dimer model. A lot of connections between dimer models, Harnack curves and amoebas are described in~\cite{Kenyon-Okounkov:2006,Kenyon-Okounkov-Sheffield:2006}

Another objects connected to the amoebas are $(p, q)$-webs~\cite{Aharony-Hanany-Kol:1998}. The spine of the amoeba of a bivariate polynomial $P$ is the $(p, q)$-web associated to the toric diagram which is the Newton polygon of $P$~\cite{Feng-He-Kennaway-Vafa:2008,Bao-He-Zahabi:2022}.

Not only amoebas but also coamoebas have their applications in dimers theory. The relationship between dimer models on the real two-torus and coamoebas of curves in $(\mathbb{C}^\times)^2$ is described in~\cite{Forsgard:2019}.

Solution to the system of fundamental thermodynamic relations in statistical thermodynamics leads to the notions of the polynomial amoeba and its contour~\cite{Passare-Pochekutov-Tsikh:2013}. It is shown in~\cite{Zahabi:2021} by using the quiver gauge theory that thermodynamic observables such as free energy, enthropy, and growth rate are explicitly derived from the limit shape of the crystal model, the boundary of the amoeba and its Harnack curve characterization. In~\cite{Konopelchenko-Angelelli:2018}, while studying the partition functions in different branches of physics, authors introduce the notion of statistical amoebas and describe their relation with polynomial amoebas.

\section{Computation of Polynomial Amoebas}

The problem of giving a complete geometric or combinatorial description for the amoeba of a polynomial has a significant computational complexity, especially for the higher dimensions. These section contains a review on the main problems of the amoeba computation and the algorithms for their solution.

\subsection{Problems and algorithms}

\subsubsection{Membership problem}

One of the basic problems of an amoeba computation is whether the given point belongs to the amoeba or, equivalently, if it belongs to a component of the amoeba complement (the membership problem). This problem is very natural, so any paper considering the computation of amoebas addresses it. Usually a solution for it is some kind of certificate~$C$ such that if $C(|z|)$ is true then Log $|z|\notin \mathcal{A}_p$ for a polynomial~$p$ (see, for example,~\cite{Forsgard-Matusevich-Mehlhop-deWolff:2018}). It is stated in~\cite[Corollary~2.7]{Theobald:2002} that for a fixed dimension~$n$ this problem can be solved in polynomial time. In~\cite{Avendano-etAl:2018} it is shown that in general this problem is {\bf PSPACE}, that is, it can be solved in polynomial time by a~parallel algorithm, provided one allows exponentially many processors.

One of the amoeba properties simplifiyng the computations is the lopsidedness (see~\cite{Forsgard:2021}). In~\cite{Purbhoo:2008} the lopsidedness criterion is presented, which provides an inequality-based certificate for non-containment of a point in an amoeba. Based on this concept, a converging sequence of approximations for the amoeba can be devised. These approximations use {\it cyclic resultants}, in~\cite{Forsgard-Matusevich-Mehlhop-deWolff:2018} a fast method of computing such resultants is provided. Another possible base for approximations like these is the real Nullstellensatz (see~\cite{Theobald-deWolff:2015}).

To solve the membership problem more efficiently, some authors propose formulations for it, that are less restrictive in some way. For example, in~\cite{Timme-Master} the soft membership problem is formulated: to determine whether the given point belongs to the amoeba with high confidence, which means that chosen criterion should only provide correct solution with some controllable high probability. This new problem is then reduced to the solution of a system of polynomial equations by means of {\it realification} of the initial polynomial (see~\cite{Theobald-deWolff:2015}). Since finding all the roots of a polynomial system has a formidable computational complexity, this problem is solved numerically by using Newton's method.

\subsubsection{Depiction of amoebas}

In two- and three-dimensional cases one of the simplest ways of describing the structure of amoeba is depicting it. For $n > 3$ it is still possible to depict three-dimensional slices of amoebas.

Since the <<tentacles>> of amoeba $\mathcal{A}_p$ always stretch to infinity, it is usually depicted only partially by choosing a domain $\Omega \subset \mathbb{R}^n$ and depicting intersection $\mathcal{A}_p \cap \Omega.$ Another complex computational problem is the choice of $\Omega$ such that the intersection inherited the essential characteristics of amoeba like the number of connected components in its complement. Algorithms for computing $\Omega$ are presented in~\cite{Timme-Master}. 

%For visualizing an amoeba we have to give the definition of the carcass of the amoeba first.

%In~\cite{Theobald2002} the homotopy-based numerical algorithm for drawing an amoeba is described. 

%\begin{definition} \rm
%{\it The carcass of the amoeba}~$\mathcal{A}$ is any subset of~$\mathcal{A}$ such that the number of connected components of the complement to the intersection~$\mathcal{A}\cap B$ for a sufficiently large ball~$B$ is as big as it could possibly be (equal to the number of connected components in~$\textrm{}^c\mathcal{A}$).
%\end{definition}
%
%The carcass of an amoeba represents its significant properties (such as the number of connected components in its complement) and may be depicted by means of graphical packages. It follows from the definition that the carcass of any given amoeba is not unique. In what follows by depicting the amoeba we mean depicting its carcass.

The basic algorithm for depicting the amoeba of a bivariate polynomial includes building a grid on~$\Omega,$ going through the grid points, replacing all the variables in the initial polynomial except for one by corresponding coordinates of the current grid point and finding roots of the obtained univariate polynomial. By substituting the roots to the grid point and applying Log map, points of amoeba are obtained. Such algorithms in different variations are presented in~\cite{Bogdanov-Kytmanov-Sadykov:2016,Forsberg-PHD,Johansson-PHD,Nilsson-PHD,Leksell-Komorowski-Bachelor}. Major drawback of this <<naive approach>> is a vast number of unnecessary computations, since not all of the obtained amoeba points belong to~$\Omega,$ and the good picture quality can require a grid with very small step because of the low control on the density of the obtained points. In~\cite{Timme-Master} another approach is presented, based on mapping of the domain~$\Omega$ to a set of pixels of an output device and executing of the membership test for corresponding points to determine whether pixel should be depicted as a part of the amoeba or not.

Some of the improvements to the computing of amoebas include construction of the grid based on Archimedean tropical hypersurfaces~\cite{Avendano-etAl:2018} which allow to execute membership tests only for the points lying in the close neighbourhood of the amoeba.  Another way to avoid testing the points which  do not belong to the amoeba is by using greedy algorithms~\cite{Timme-Master}, that is, testing only points in some neighborhood of already tested ones. 

In~\cite{Zhukov-Sadykov:2023} authors compute the connected components of the amoeba complement by using their order. Since its value is the same for all points belonging to the same component and differs for points from different components, the order of the component is a good classifier for points of an amoeba complement. For points of the amoeba the integral for computing the order does not always converge but the jump of the order function itself can be a criterion of a point belonging to the amoeba.

In~\cite{Bao-He-Hirst:2023} computational problems for amoebas are solved by the means of machine learning techniques. Authors use artificial neural networks to determine the genus of an amoeba and to solve the membership problem. There are examples in this article presenting models with prediction accuracy around 0,95.

\subsection{Software}

Let us consider three of the latest software examples for depicting the amoebas and compare their functionality. All of these solutions are freeware and available online. 

{\bf [a]} The script by Dmitry Bogdanov at \url{http://dvbogdanov.ru/amoeba}. Very simple, both in usage and its functionality, it allows to depict amoebas and compactified amoebas for input polynomials. The script itself only generates the MatLab code for depicting the amoebas with given parameters. This approach has the disadvantage of needing a Matlab installation, but the benefit of being hardware-independent. Examples of the most complex amoebas depicted with this software are presented at~\url{www.researchgate.net/publication/338341129_Giant_amoeba_zoo}.

{\bf [b]} Package PolynomialAmoebas written by Sascha Timme in the Julia programming language (\url{https://github.com/saschatimme/PolynomialAmoebas.jl}. It has very broad functionality, including not only tools for depicting the amoebas and not only in $2$ dimensions. It also allows one to depict spines and contours of amoebas, coamoebas and amoebas in~$3$ dimensions.   

{\bf [c]} The project by Timur Sadykov and Timur Zhukov (\url{http://amoebas.ru/index.html}). It includes the visualization of amoebas in $2$ and $3$ dimensions, three-dimensional slices of 4D amoebas and allows to visualize the evolution of amoeba due to a change of the coefficients in the polynomial. The project is still under development.

\subsection{Software comparison}

%In Table~\ref{tab:CompareFunctionality} we present comparison of the functionality of the software tools mentioned above.

%In Table 1 we present comparison of the functionality of the software tools mentioned above.
%\vskip0.2cm
%%\begin{table}[h!]
%%\caption{Comparison of the functionality of the packages for computing the polynomial amoebas}
%\noindent Table 1: Comparison of the functionality of the packages for computing the polynomial amoebas
%
%\noindent\begin{tabular}{|p{6cm}|wc{2cm}|wc{2cm}|wc{2cm}|}
%\hline
%\centering Functionality & [a] & [b] & [c] \\
%\hline
%Newton polytopes & + & & + \\ 
%\hline
%Amoebas &  + &   + &  + \\
%\hline
%Compactified amoebas & + & & \\
%\hline
%Contours of amoebas & & + & \\
%\hline 
%Spines of amoebas & & + & \\
%%\hline
%%Coamoebas & & + & \\
%\hline
%3D amoebas & & + & + \\
%\hline
%Amoeba evolution & & & + \\
%\hline
%\end{tabular}
%\vskip0.2cm
%\label{tab:CompareFunctionality}
%\end{table}

%In Table~\ref{tab:ComparePerformance} the performance of these packages is compared. Seconds are chosen to be the units of measure there on purpose, since the time in milliseconds may vary slightly in different runs.

The hardware configuration used for benchmarking is as follows: Intel Xeon Gold 6146 CPU with a clock speed of 3.2 GHz and 128 GB RAM. To measure the time required to run a function for the software [a] MatLab IDE built-in functionality was used, for the software [b] -- @btime function from BenchmarkTools package and for the software [c] -- built-in functionality in the desktop version.

The packages differ from each other significantly, both in the implementation and the functionality, so the important question for measuring the performance is a choice of parameters. For example, depicting amoebas with software [a] implies a choice of grid by choosing numbers of modulus values~$r_x, r_y$ and a number of argument values~$\phi.$ Detailed information on the algorithm and its parameters is given in~\cite{Bogdanov-Kytmanov-Sadykov:2016}. For the low values of these parameters, algorithm terminates faster, but the picture of amoeba dissipates into the set of isolated points, especially after zooming in it. For the tests, 2 sets of the parameters are used -- the ones by default (each number equals 100) and when each number equals 500, which results in smooth picture. For the simplicity in what follows these cases are denoted by <<100 values>> and <<500 values>> respectively.

For software [c] the important parameter affecting the time of termination and the quality of resulting picture is the number of iterations the algorithm performs. %In Table~\ref{tab:ComparePerformance} for each polynomial results for two values of this parameter are presented, first of them is a minimal one so that the number of connected components in the result coincides with actual number for the amoeba, second is the value for which the picture looks smooth enough. %As one may see in the table, the amoeba for the last polynomial is only successfully depicted by the first package, the reason for that is probably the issue of working with big numbers, which MatLab obviously avoids. 

\begin{example} \rm
Consider the polynomial~$p_2(z_1,z_2) = 5 z_1 + 15 z_1^2 + 240 z_1 z_2 + 400 z_1^2 z_2 + 8 z_1^3 z_2 + 10 z_2^2 + 900 z_1 z_2^2 + 900 z_1^2 z_2^2 + 10 z_1^3 z_2^2 + 8 z_2^3 + 400 z_1 z_2^3 + 240 z_1^2 z_2^3 + 15 z_1 z_2^4 + 5 z_1^2 z_2^4.$ Let us compute the amoeba for it using the software [a], [b] and [c].

The result obtained with the software [a] for 100 values is shown in Figure~\ref{fig:p2_software_a} (a). It describes the structure of the amoeba, but the quality of picture is not high enough. The computation time equals~1.6 seconds. For 500 values the computation time becomes 16 seconds, and there are no points in the picture that look isolated (see Figure~\ref{fig:p2_software_a} (b)). 

\begin{figure}[h!]
\hskip0.5cm
\begin{minipage}{5cm}
\includegraphics[width=5cm]{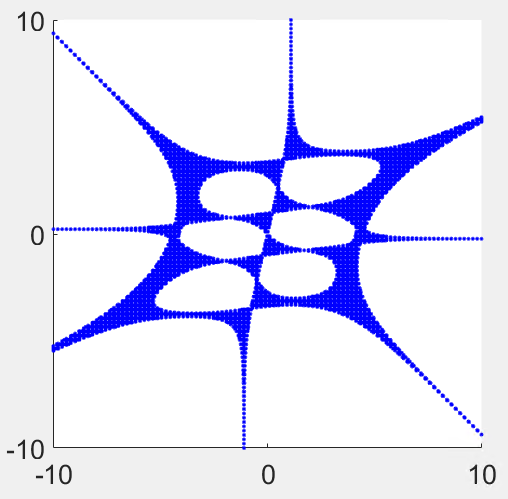}
\centering $(a)$
\end{minipage}
\hskip1cm
\begin{minipage}{5cm}
\includegraphics[width=5cm]{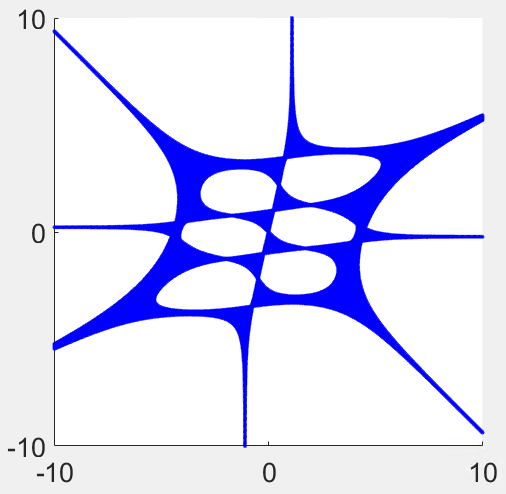}
\centering $(b)$
\end{minipage}

\caption{The amoeba for the polynomial~$p_2$ computed with the software [a] -- (a) 100 values, (b) 500 values}
\label{fig:p2_software_a}
\end{figure}

The computation result for the software [b] with default parameters is presented in Figure~\ref{fig:p2_software_b}. This package also provides faster algorithms (for example, the greedy one), but the picture for the default algorithm has the best quality.

\begin{figure}[h!]
\hskip0.5cm
\begin{minipage}{6cm}
\includegraphics[width=6cm]{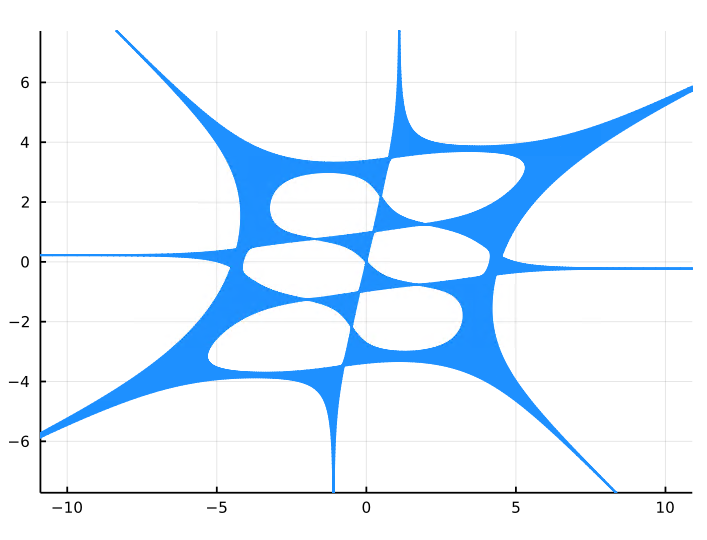}
\end{minipage}

\caption{The amoeba for the polynomial~$p_2$ computed with the software [b]}
\label{fig:p2_software_b}
\end{figure}

The result for the software [c] is shown in Figure~\ref{fig:p2_software_c} for 2 different values of iterations parameter: 5 cycles, which is the minimal value such that the number of connected components of the amoeba complement in the picture coincides with the actual number (computation terminates in 5 seconds), and 9 cycles, when the picture is smooth (computation terminates in 450 seconds).

\begin{figure}[h!]
\hskip0.5cm
\begin{minipage}{5cm}
\includegraphics[width=5cm]{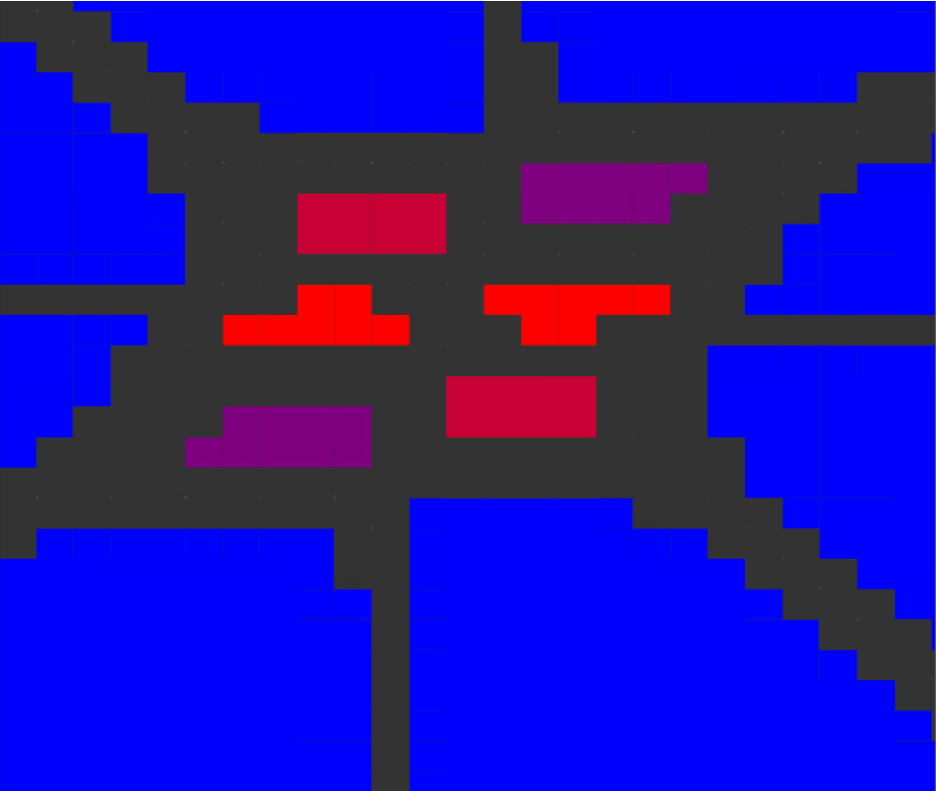}
\centering $(a)$
\end{minipage}
\hskip1cm
\begin{minipage}{5cm}
\includegraphics[width=5cm]{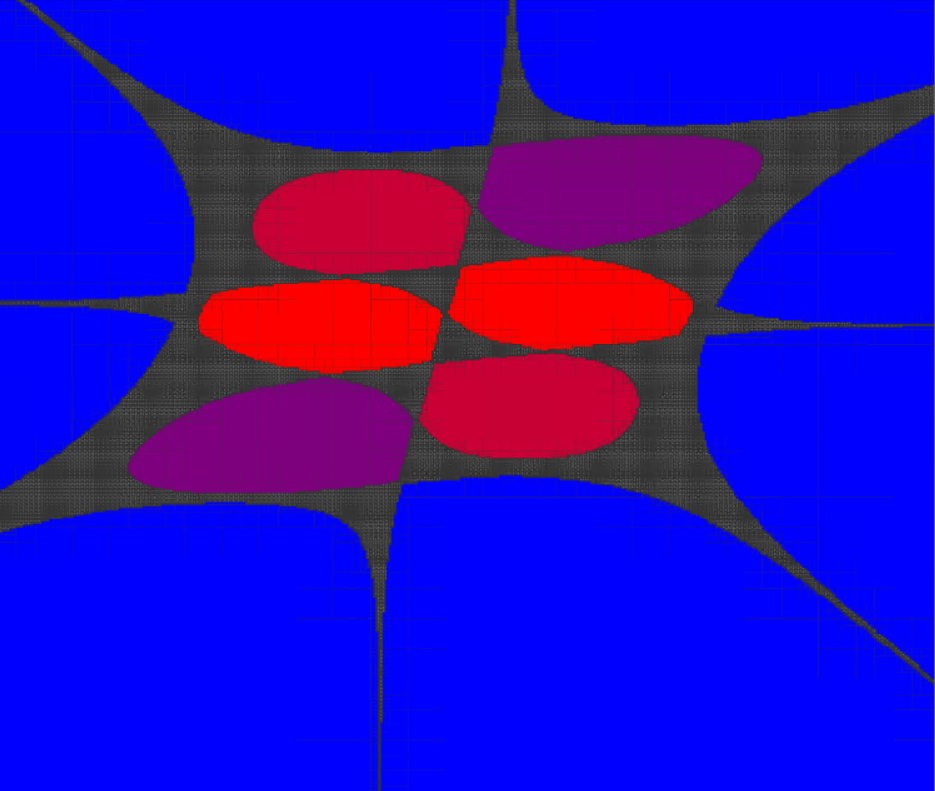}
\centering $(b)$
\end{minipage}

\caption{The amoeba for the polynomial~$p_2$ computed with the software [c] -- (a) 5 cycles; (b) 9 cycles}
\label{fig:p2_software_c}
\end{figure}
\end{example}

\subsubsection{Tests}

\begin{figure}[h!]
\begin{minipage}{6.5cm}
\includegraphics[width=7cm]{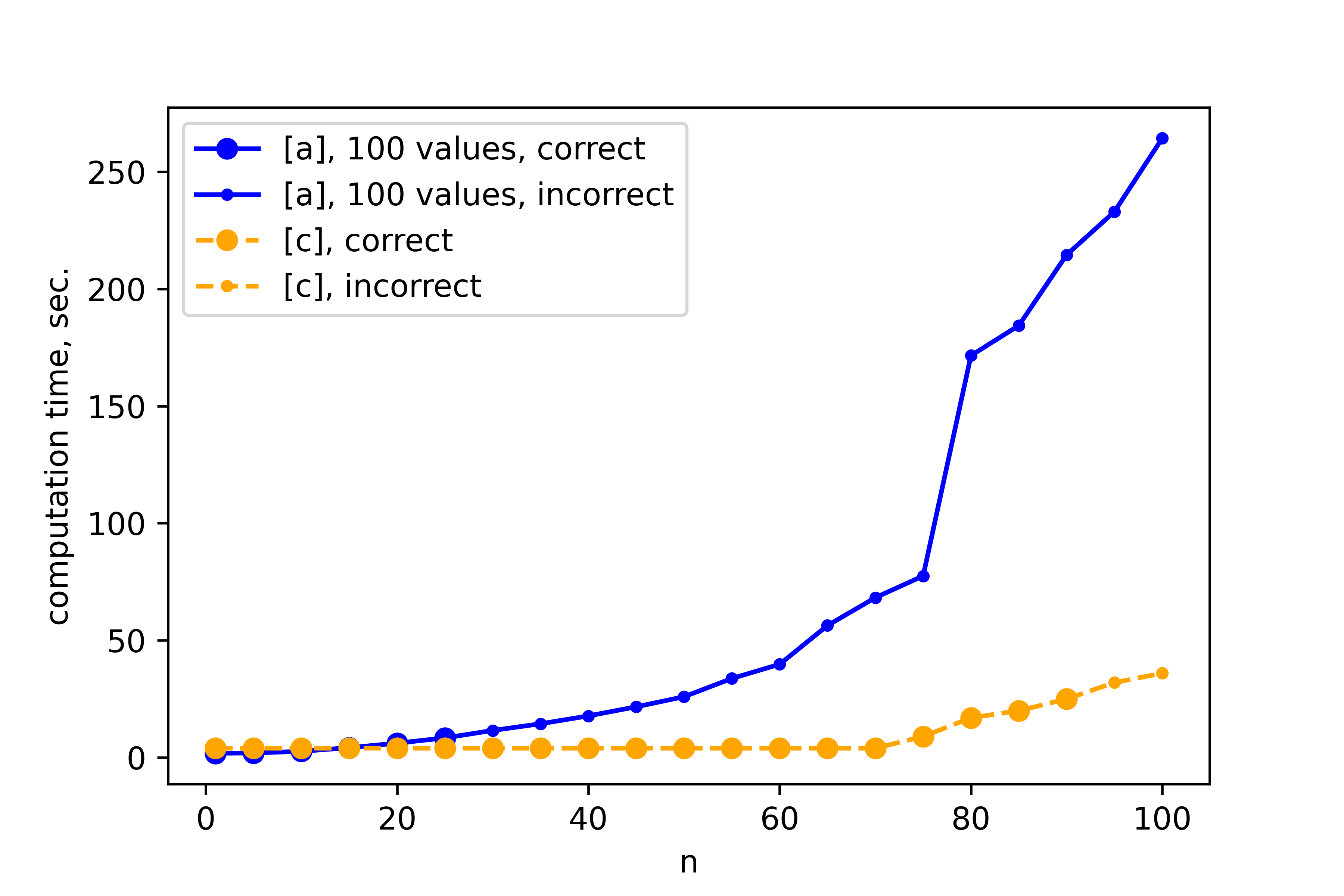}
\centering $(a)$
\end{minipage}
\begin{minipage}{6.5cm}
\includegraphics[width=7cm]{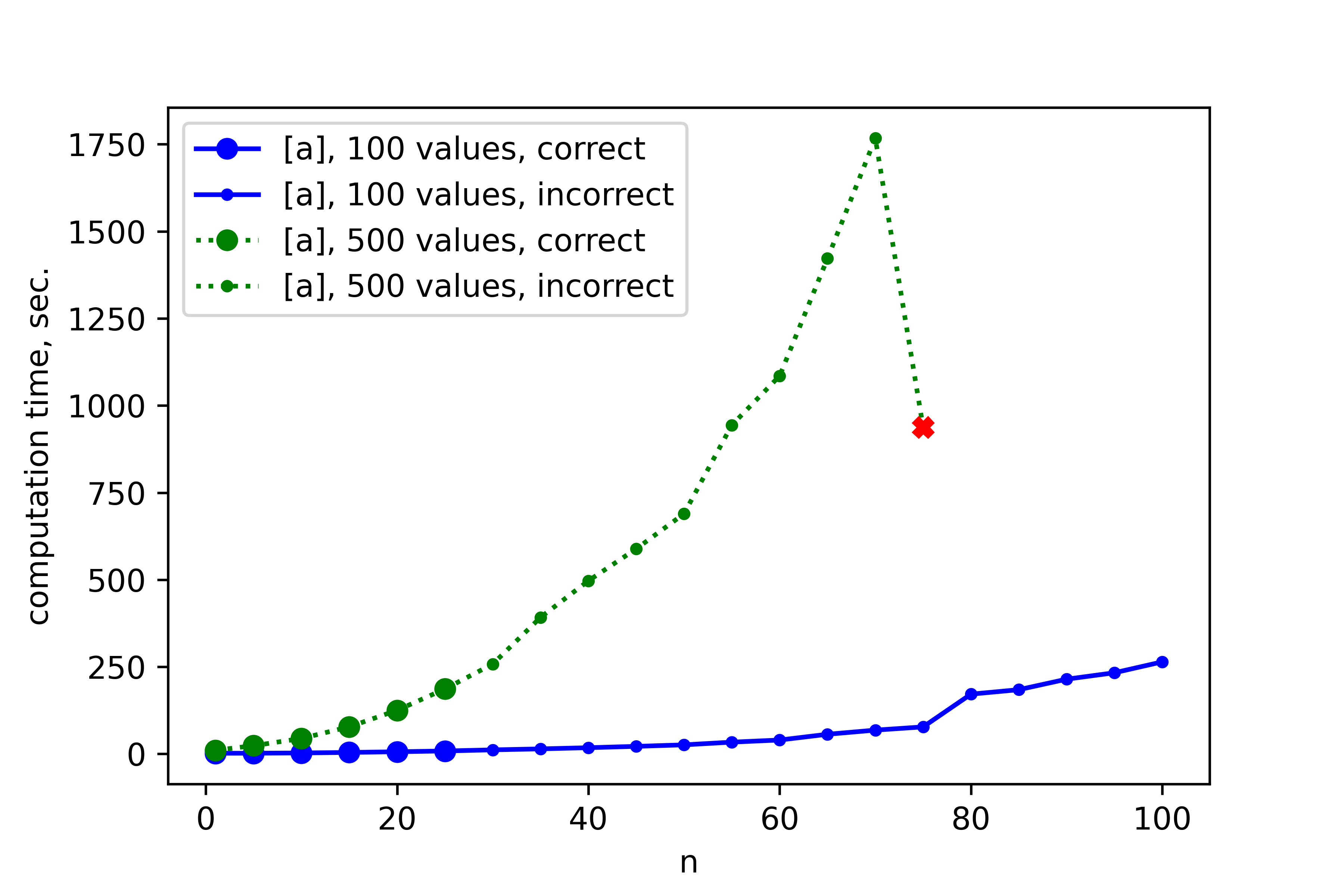}
\centering $(b)$
\end{minipage}

\caption{Results of the test computations on polynomials $z_1^n+z_2^n+1$ -- (a)~comparison of the software [a] for 100 values and the software [c]; (b)  comparison of the software [a] for 100 values and for 500 values}
\label{fig:test_1}
\end{figure}

First test is performed on the polynomials of the form $z_1^n+z_2^n+1, n=1,2,\ldots$ to determine a dependence of the performance of the software on the degree of the initial polynomial. These polynomials are maximally sparse and their amoebas are solid. Computation results are presented in Figure~\ref{fig:test_1}. Even for trinomials, starting with some value of~$n$ amoeba pictures contain image artifacts, these cases are represented by smaller points in the graph. In what follows we refer to these cases as incorrect ones.

The comparison of the performance of the software [a] and the software [c] is presented in Figure~\ref{fig:test_1} (a). For the software [a] image becomes incorrect starting with $n=30,$ for the software [c] -- with $n=95.$ In general, the software [c] shows better results in this test -- it is faster and performs correct computations for larger number of~$n$ values.

In Figure~\ref{fig:test_1} (b) performance of the software [a] for the cases of 100 values and 500 values are compared. Image artifacts in both cases first appear for~$n=30.$ Starting from $n=75$ computations for 500 values terminate with an error.

The software [b] seems to be inappropriate for the mass tests, since a lot of polynomials given as an input just lead to generating an exception. In particular, this applies to the polynomials $z_1^n+z_2^n+1.$ Possible reason for this is an update of the Julia packages combined with the lack of a support from the developer.

\begin{figure}[h!]
\begin{minipage}{6.5cm}
\includegraphics[width=7cm]{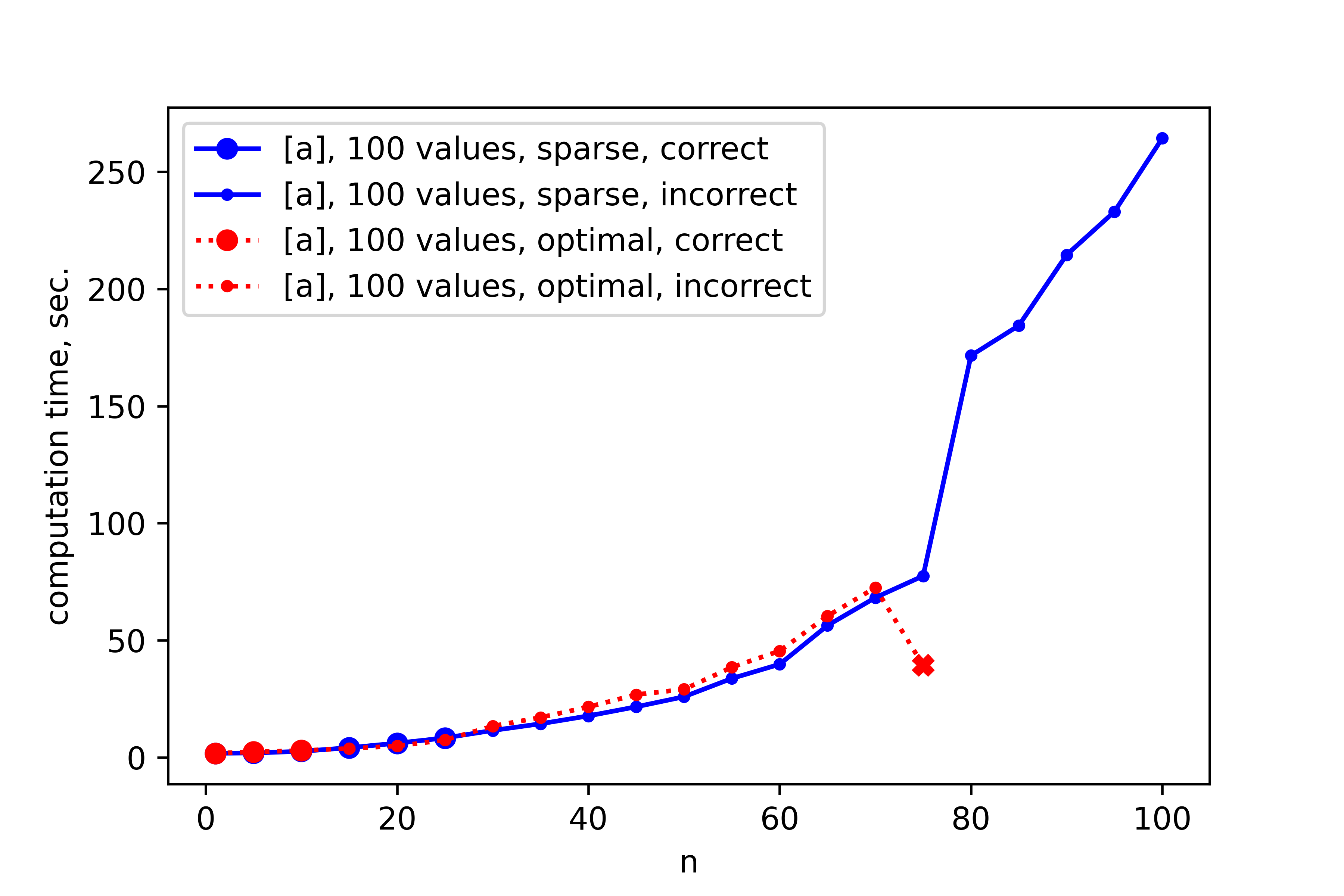}
\centering $(a)$
\end{minipage}
\begin{minipage}{6.5cm}
\includegraphics[width=7cm]{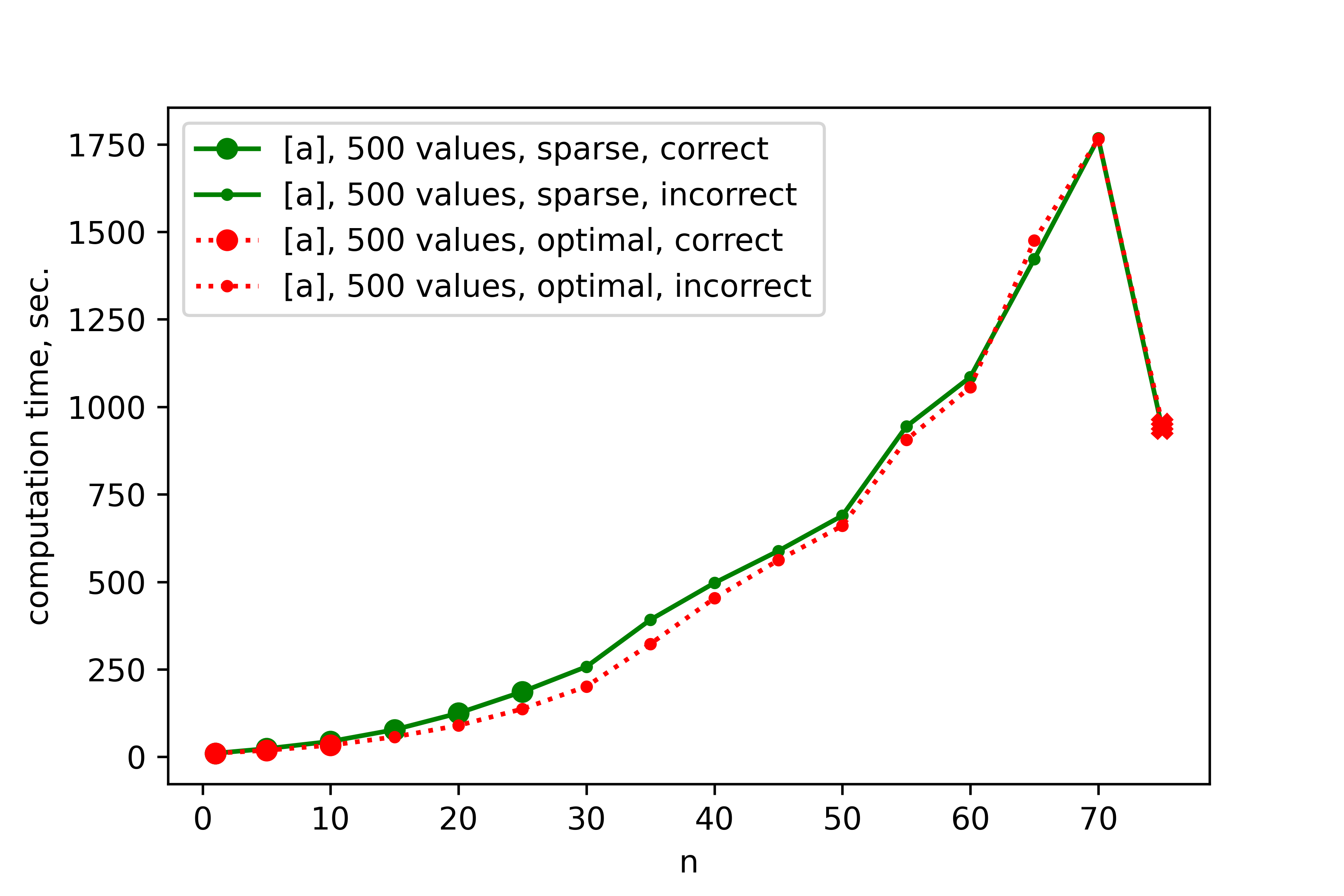}
\centering $(b)$
\end{minipage}
\begin{minipage}{6.5cm}
\includegraphics[width=7cm]{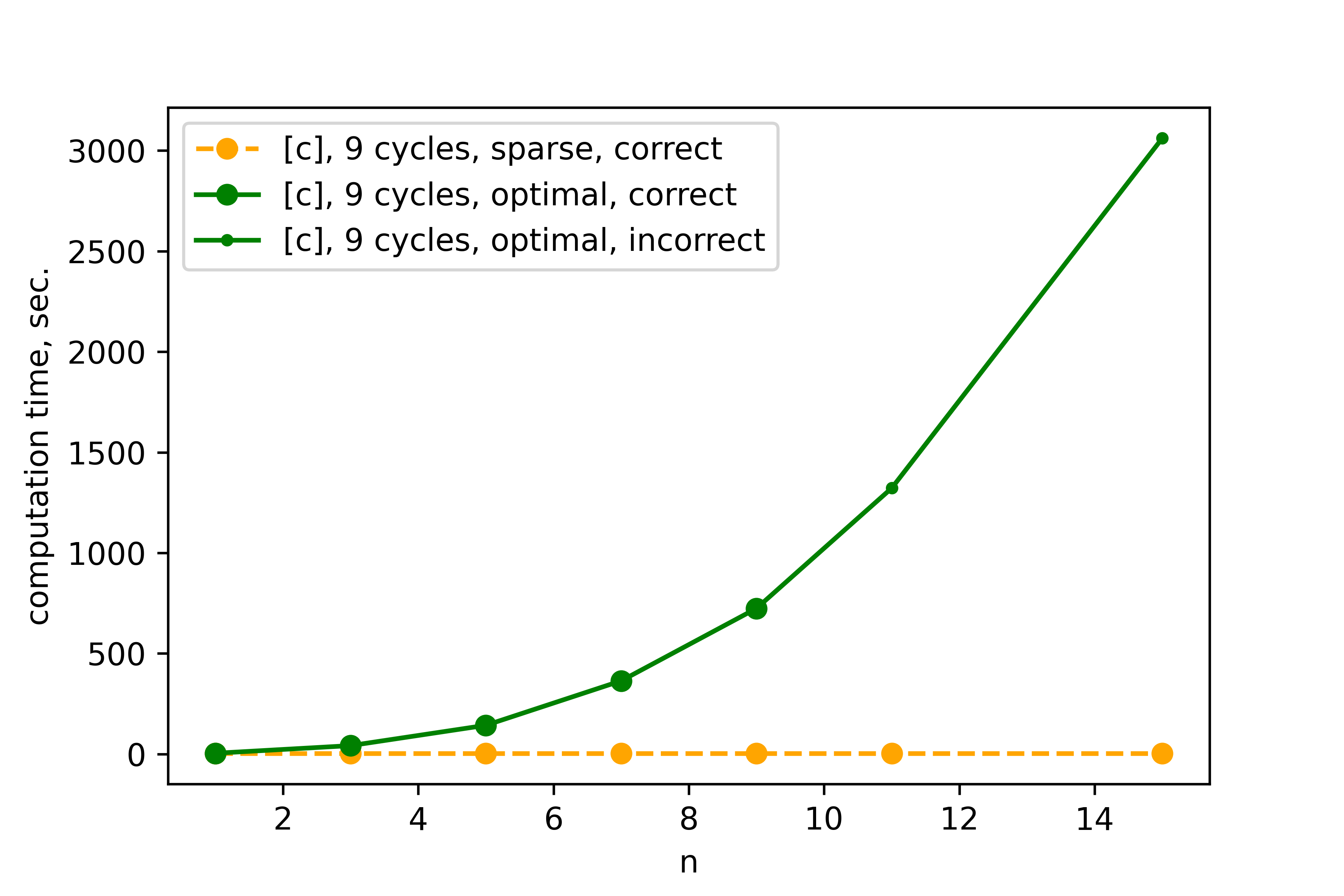}
\centering $(c)$
\end{minipage}
\caption{Results of the test computations on optimal polynomials -- (a) comparison of the software [a] for 100 values on maximally sparse and optimal polynomials; (b) comparison of the software [a] for 500 values on maximally sparse and optimal polynomials; (c) comparison of the software [c] for 9 cycles on maximally sparse and optimal polynomials}
\label{fig:test_2}
\end{figure}

Second test is performed on the optimal polynomials supported in simplexes, generated by the means of an algorithm presented in~\cite{Zhukov-Sadykov:2023}. Results are presented in Figure~\ref{fig:test_2}, where (a) and (b) show the comparison of computation times for the software [a] for sparse polynomials from the previous test and optimal polynomials. Surprisingly, this time does not depend on the number of monomials at all, in the case of 500 values the computation time for optimal polynomials for some interval of values of $n$ is even lesser than for maximally sparse ones. The only difference is that the algorithm does not terminate with an error for $n=75$ in the case of 100 values.

It must be noted that for optimal polynomials smaller points on the graph denote the computations such that the algorithm do not recognize all of the connected components of an amoeba complement, since these components are too small. Maximal value of~$n$ such that the algorithm recognizes all of the components is $n=10$ in the case of 100 values and $n=11$ in the case of 500 values.

Figure~\ref{fig:test_2} (c) depicts the results for the software [c] -- and there the computation time grows fast with the growth of~$n.$ The algorithm for 9 cycles fails to recognize all of the components of the amoeba complement starting at~$n=11.$

\begin{figure}[h!]
\begin{minipage}{6.5cm}
\includegraphics[width=7cm]{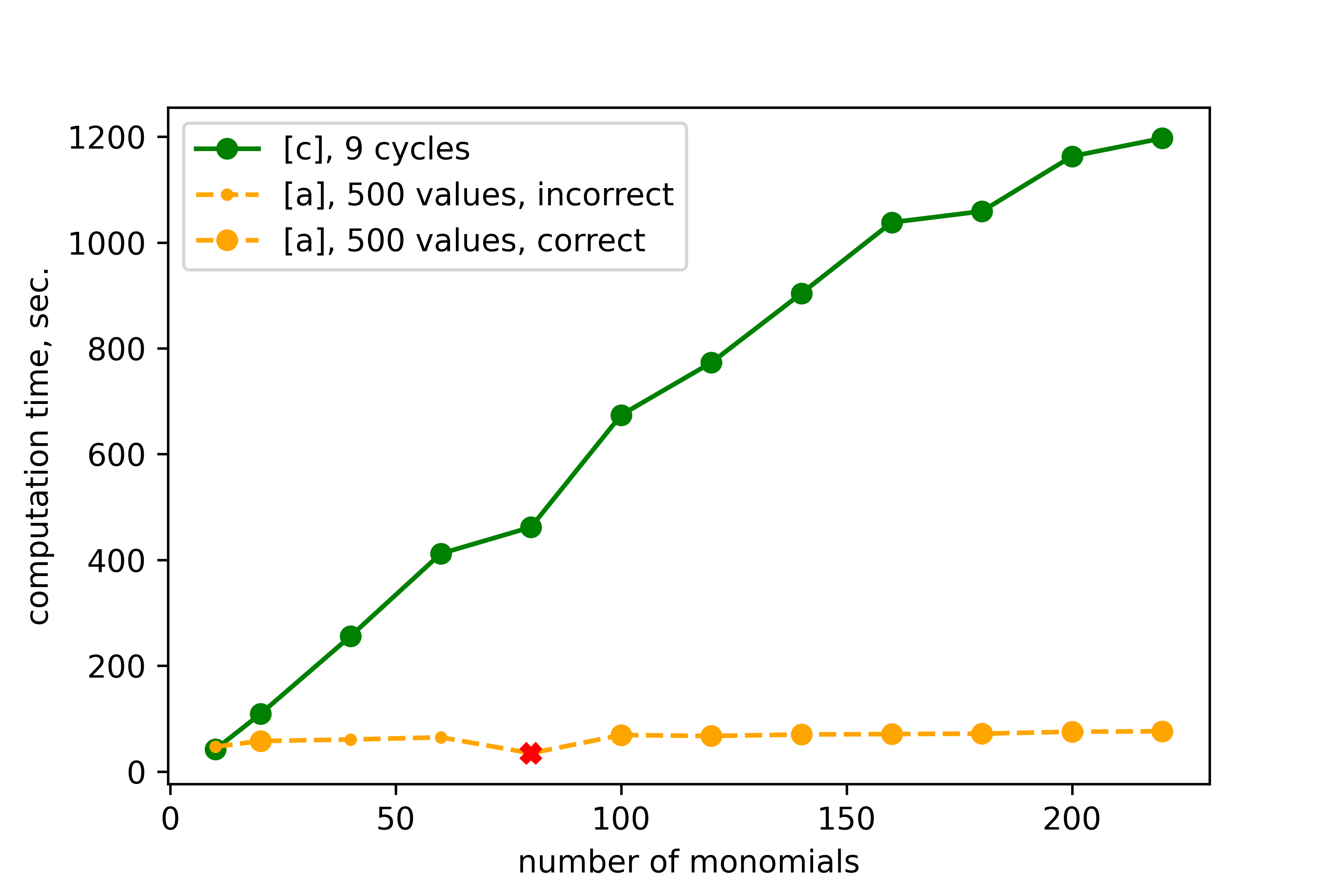}
\end{minipage}
\caption{Results of the test computations on random polynomials -- comparison of the software [a] and the software [c]}
\label{fig:test_3}
\end{figure}

The final test was performed on the sequence of polynomials with random integer coefficients in $[1,5]$ and the increasing number of monomials. The computation results for the software [a] and the software [c] are presented in Figure~\ref{fig:test_3}. All of the graph points for the software [c] are depicted large, since in this case it is not easy to understand whether the algorithm recognizes all of the components of the amoeba complement. The randomness of input polynomials leads to highly unpleasant amoebas with a lot of tentacles and small spaces between them. Results with the software [a] in this case are also random, for some monomial numbers (especially, the lower ones) the resulting amoeba contained obvious artifacts and for 80 monomials algorithm terminated with an error.

\subsubsection{Test results}

The software [b] in the present state fits only for computing in some particular cases, since some input polynomials only generate an error and there were some cases, when the algorithm went into the infinite loop for no obvious reasons. For suitable polynomials the computation time is better than for the software [c].

The software [a] is faster than the software [c] in the case of optimal polynomials, but it has limitations on the degree of a polynomial, which becomes stricter with a growth of the number of monomials. For random polynomials there are a lot of issues (including image artifacts and even computation errors) with the result. 

The software [c] for the parameter values which ensure the best performance is slower than the software [a] and the computation time depends on the number of monomials, but it has lesser limitations on the parameters of input polynomials than other two packages.

 \bibliographystyle{splncs04}
 \bibliography{amoebas}

\begin{thebibliography}{10}
\providecommand{\url}[1]{\texttt{#1}}
\providecommand{\urlprefix}{URL }
\providecommand{\doi}[1]{https://doi.org/#1}

\bibitem{PassareConjectureFormulation:2009}
Open problems: Amoebas and tropical geometry (2009),
  \url{http://aimath.org/mathresources/openproblems/}

\bibitem{Aharony-Hanany-Kol:1998}
Aharony, O., Hanany, A., Kol, B.: Webs of $(p, q)$ 5-branes, five dimensional
  field theories and grid diagrams. JHEP  \textbf{1} (1998)

\bibitem{Avendano-etAl:2018}
Avenda$\tilde{\textrm{n}}$o, M., Kogan, R., Nisse, M., Rojas, J.: Metric
  estimates and membership complexity for archimedean amoebae and tropical
  hypersurfaces. Reviews in Mathematical Physics  \textbf{46},  45--65 (2018)

\bibitem{Bao-He-Zahabi:2022}
Bao, J., He, Y.H., Zahabi, A.: Mahler measure for a quiver symphony.
  Communications in Mathematical Physics  (2022)

\bibitem{Bao-He-Hirst:2023}
Bao, J., He, Y.H., Hirst, E.: Neurons on amoebae. Journal of Symbolic
  Computation  \textbf{116},  1--38 (2023)

\bibitem{Bogdanov-Kytmanov-Sadykov:2016}
Bogdanov, D., Kytmanov, A., Sadykov, T.: Algorithmic computation of polynomial
  amoebas. Lecture Notes in Computer Science (including subseries Lecture Notes
  in Artificial Intelligence and Lecture Notes in Bioinformatics  \textbf{9890
  LNCS},  87--100 (2016)

\bibitem{Bogdanov-Sadykov:2020}
Bogdanov, D., Sadykov, T.: Hypergeometric polynomials are optimal.
  Mathematische Zeitschrift  \textbf{296}(1-2),  373--390 (2020)

\bibitem{Cauchy-OeuvresCompletes}
Cauchy, O.: {\OE}uvres compl\'etes d'Augustin Cauchy. Addison Wesley,
  Massachusetts, 2 edn. (1882--1938)

\bibitem{Cherepanskiy-Tsikh:2020}
Cherepanskiy, A., Tsikh, A.: Convergence of two-dimensional hypergeometric
  series for algebraic functions. Integral Transforms and Special Functions
  \textbf{31}(10),  838--855 (2020)

\bibitem{Draisma-Rau-Yuen:2020}
Draisma, J., Rau, J., Yuen, C.: The dimension of an amoeba. Bulletin of the
  London Mathematical Society  \textbf{52}(1),  16--23 (2020)

\bibitem{Einsiedler-Kapranov-Lind:2006}
Einsiedler, M., Kapranov, M., Lind, D.: Non-archimedean amoebas and tropical
  varieties. Journal fur die Reine und Angewandte Mathematik  \textbf{601},
  139--157 (2006)

\bibitem{Feng-He-Kennaway-Vafa:2008}
Feng, B., He, Y.H., Kennaway, K., Vafa, C.: Dimer models from mirror symmetry
  and quivering amoeb\ae. Advances in Theoretical and Mathematical Physics
  \textbf{12}(3),  489--545 (2008)

\bibitem{Forsberg-PHD}
Forsberg, M.: Amoebas and Laurent Series. Ph.D. thesis, Royal Institute of
  Technology, Stockholm, Sweden (1998)

\bibitem{Forsberg-Passare-Tsikh:2000}
Forsberg, M., Passare, M., Tsikh, A.: Laurent determinants and arrangements of
  hyperplane amoebas. Advances in Mathematics  \textbf{151},  45--70 (2000)

\bibitem{Forsgard:2019}
Forsg{\aa}rd, J.: On dimer models and coamoebas. Annales de l'Institut Henri
  Poincare (D) Combinatorics, Physics and their Interactions  \textbf{6}(2),
  199--219 (2019)

\bibitem{Forsgard:2020}
Forsg{\aa}rd, J.: Tropical approximation of exponential sums and the
  multivariate $\textrm{F}$ujiwara bound. Moscow Mathematical Journal
  \textbf{20}(2),  311--321 (2020)

\bibitem{Forsgard:2021}
Forsg{\aa}rd, J.: Discriminant amoebas and lopsidedness. Journal of Commutative
  Algebra  \textbf{13}(1),  41--60 (2021)

\bibitem{Forsgard-Matusevich-Mehlhop-deWolff:2018}
Forsg{\aa}rd, J., Matusevich, L., Mehlhop, N., de~Wolff, T.: Lopsided
  approximation of amoebas. Mathematics of Computation  \textbf{88}(315),
  485--500 (2018)

\bibitem{GKZ:1994}
Gelfand, I., Kapranov, M., Zelevinsky, A.: Discriminants, resultants, and
  multidimensional determinants. Birkh\"auser Boston Inc., Boston, MA (1994)

\bibitem{Goucha-Gouveia:2021}
Goucha, A., Gouveia, J.: The phaseless rank of a matrix. SIAM Journal on
  Applied Algebra and Geometry  \textbf{5}(3),  526--551 (2021)

\bibitem{Guilloux-Marche:2021}
Guilloux, A., March\'e, J.: Volume function and $\textrm{M}$ahler measure of
  exact polynomials. Compositio Mathematica pp. 809--834 (2021)

\bibitem{Harnack:1876}
Harnack, A.: $\textrm{Ü}$ber vieltheiligkeit der ebenen algebraischen curven.
  Math. Ann.  \textbf{10},  189--199 (1876)

\bibitem{Hicks:2020}
Hicks, J.: Tropical lagrangian hypersurfaces are unobstructed. Journal of
  Topology  \textbf{13}(4),  1409--1454 (2020)

\bibitem{Hilbert:1902}
Hilbert, D.: Mathematical problems. Bull. Amer. Math. Soc.  \textbf{8},
  437--479 (1902)

\bibitem{Iliman-deWolff:2016}
Iliman, S., de~Wolff, T.: Amoebas, nonnegative polynomials and sums of squares
  supported on circuits. Research in Mathematical Sciences  \textbf{3}(1)
  (2016)

\bibitem{Jensen-Leykin-Yu:2016}
Jensen, A., Leykin, A., Yu, J.: Computing tropical curves via homotopy
  continuation. Experimental Mathematics  \textbf{25}(1),  83--93 (2016)

\bibitem{Johansson-PHD}
Johansson, P.: On the topology of the coamoeba. Ph.D. thesis, Stockholm
  University, Sweden (2014)

\bibitem{Johansson-Samuelsson:2017}
Johansson, P., Kalm, H.S.: A $\textrm{R}$onkin type function for coamoebas.
  Journal of Geometric Analysis  \textbf{27}(1),  643--670 (2017)

\bibitem{Jonsson:2016}
Jonsson, M.: Degenerations of amoebae and $\textrm{B}$erkovich spaces.
  Mathematische Annalen  \textbf{364}(1-2),  293--311 (2016)

\bibitem{Juhnke-Kubitzke-deWolff:2016}
Juhnke-Kubitzke, M., de~Wolff, T.: Intersections of amoebas. In: 28th
  International Conference on Formal Power Series and Algebraic Combinatorics,
  FPSAC 2016. pp. 659--670. Vancouver, USA (2016)

\bibitem{Kaveh-Manon:2019}
Kaveh, K., Manon, C.: Gr\"obner theory and tropical geometry on spherical
  varieties. Transformation Groups  \textbf{24}(4),  1095--1145 (2019)

\bibitem{Kenyon:2008}
Kenyon, R.: An introduction to the dimer model (2008),
  \url{https://arxiv.org/pdf/math/0310326.pdf}

\bibitem{Kenyon-Okounkov:2006}
Kenyon, R., Okounkov, A.: Planar dimers and $\textrm{H}$arnack curves. Duke
  Mathematical Journal  \textbf{131}(3),  499--524 (2006)

\bibitem{Kenyon-Okounkov-Sheffield:2006}
Kenyon, R., Okounkov, A., Sheffield, S.: Dimers and amoebae. Annals of
  Mathematics  \textbf{163}(3),  1019--1056 (2006)

\bibitem{Khovanskii:1991}
Khovanskii, A.: Translations of Mathematical Monographs. Volume 88. Fewnomials.
  American Mathematical Society (1991)

\bibitem{Kim-Nisse:2021}
Kim, Y., Nisse, M.: A natural topological manifold structure of phase tropical
  hypersurfaces. Journal of the Korean Mathematical Society  \textbf{58}(2),
  451--471 (2021)

\bibitem{Konopelchenko-Angelelli:2018}
Konopelchenko, B., Angelelli, M.: Zeros and amoebas of partition functions.
  Reviews in Mathematical Physics  \textbf{30}(9) (2018)

\bibitem{Lang:2019}
Lang, L.: Amoebas of curves and the $\textrm{L}$yashko–$\textrm{L}$ooijenga
  map. Journal of the London Mathematical Society  \textbf{100}(1),  301--322
  (2019)

\bibitem{Lang:2020}
Lang, L.: Harmonic tropical morphisms and approximation. Mathematische Annalen
  \textbf{377}(1-2),  379--419 (2020)

\bibitem{Lang-Shapiro-Shustin:2021}
Lang, L., Shapiro, B., Shustin, E.: On the number of intersection points of the
  contour of an amoeba with a line. Indiana University Mathematics Journal
  \textbf{70}(4),  1335--1353 (2021)

\bibitem{Leksell-Komorowski-Bachelor}
Leksell, M., Komorowski, W.: Amoeba program: Computing and visualizing amoebas
  for some complex-valued bivariate expressions (2007)

\bibitem{Lyapin:2009}
Lyapin, A.: Riordan's arrays and two-dimensional difference equations. Journal
  of Siberian Federal University. Mathematics and Physics  \textbf{2}(2),
  210--220 (2009)

\bibitem{Mikhalkin:2000}
Mikhalkin, G.: Real algebraic curves, the moment map and amoebas. Annals of
  Mathematics pp. 309--326 (2000)

\bibitem{Mikhalkin:2004}
Mikhalkin, G.: Amoebas of algebraic varieties and tropical geometry. In:
  Different faces of geometry. Int. Math. Series (N.Y) 3, pp. 257--300. Kluwer
  (2004)

\bibitem{Mikhalkin:2006}
Mikhalkin, G.: Tropical geometry and its applications. In: International
  Congress of Mathematicians, ICM 2006. vol.~2, pp. 827--852 (2006)

\bibitem{Mikhalkin:2017}
Mikhalkin, G.: Amoebas of half-dimensional varieties. Trends in Mathematics
  \textbf{9783319524696},  349--359 (2017)

\bibitem{Mikhalkin-Rullgard:2001}
Mikhalkin, G., Rullg{\aa}rd, H.: Amoebas of maximal area. Internat. Math. Res.
  Notices pp. 441--451 (2001)

\bibitem{Mikhalkin-Shkolnikov:2022}
Mikhalkin, G., Rullg{\aa}rd, H.: Non-commutative amoebas. Bulletin of the
  London Mathematical Society  \textbf{54}(2),  335--368 (2022)

\bibitem{Mkrtchian-Yuzhakov:1982}
Mkrtchian, M., Yuzhakov, A.: The $\textrm{N}$ewton polytope and the
  $\textrm{L}$aurent series of rational functions of $n$ variables. Izv. Akad.
  Nauk ArmSSR  \textbf{17},  99--105 (1982)

\bibitem{Nilsson-PHD}
Nilsson, L.: Amoebas, discriminants and hypergeometric functions. Ph.D. thesis,
  Stockholm University, Sweden (2009)

\bibitem{Nilsson-Passare-Tsikh:2019}
Nilsson, L., Passare, M., Tsikh, A.: Domains of convergence for ahypergeometric
  series and integrals. Journal of Siberian Federal University - Mathematics
  and Physics  \textbf{12}(4),  509--529 (2019)

\bibitem{Nisse:2016}
Nisse, M.: Amoeba basis of zero-dimensional varieties. Journal of Pure and
  Applied Algebra  \textbf{220}(3),  1252--1257 (2016)

\bibitem{Nisse-Sadykov:2019}
Nisse, M., Sadykov, T.: Amoeba-shaped polyhedral complex of an algebraic
  hypersurface. Journal of Geometric Analysis  \textbf{29}(2),  1356--1368
  (2019)

\bibitem{Passare:2016}
Passare, M.: The trigonometry of $\textrm{H}$arnack curves. Journal of Siberian
  Federal University - Mathematics and Physics  \textbf{9}(3),  347--352 (2016)

\bibitem{Passare-Pochekutov-Tsikh:2013}
Passare, M., Pochekutov, D., Tsikh, A.: Amoebas of complex hypersurfaces in
  statistical thermodynamics. Mathematical Physics, Analysis and Geometry
  \textbf{16},  89--108 (2013)

\bibitem{Passare-Rullgard:2000}
Passare, M., Rullgård, H.: Amoebas, $\textrm{M}$onge-$\textrm{A}$mp\'ere
  measures and triangulations of the $\textrm{N}$ewton polytope. Duke
  Mathematical Journal  \textbf{121}(3),  481--507 (2004)

\bibitem{Passare-Tsikh:2002}
Passare, M., Sadykov, T., Tsikh, A.: Algebraic equations and hypergeometric
  series. In: The Legacy of Niels Henrik Abel, pp. 653--672. The Abel
  Bicentennial, Oslo (2002)

\bibitem{Passare-Sadykov-Tsikh:2005}
Passare, M., Sadykov, T., Tsikh, A.: Singularities of hypergeometric functions
  in several variables. Compositio Mathematica  \textbf{141}(3),  787--810
  (2005)

\bibitem{Purbhoo:2008}
Purbhoo, K.: A $\textrm{N}$ullstellensatz for amoebas. Duke Math. J.
  \textbf{141}(3),  407--445 (2008)

\bibitem{Ruan:2000}
Ruan, W.D.: Newton polygon and string diagram (2000)

\bibitem{Theobald:2002}
Theobald, T.: Computing amoebas. Experiment. Math.  \textbf{11}(4),  513--526
  (2002)

\bibitem{Theobald-deWolff:2015}
Theobald, T., de~Wolff, T.: Approximating amoebas and coamoebas by sums of
  squares. Math. of Computation  \textbf{84}(291),  455--473 (2015)

\bibitem{Theobald-deWolff:2016}
Theobald, T., de~Wolff, T.: Norms of roots of trinomials. Mathematische Annalen
   \textbf{366}(1-2),  219--247 (2016)

\bibitem{Timme-Master}
Timme, S.: Fast computation of amoebas, coamoebas and imaginary projections in
  low dimensions (2018)

\bibitem{deWolff:2017}
de~Wolff, T.: Amoebas and their tropicalizations -- a survey. In: Analysis
  Meets Geometry, pp. 157--190 (2017)

\bibitem{Xu-Dunkl:2014}
Xu, Y., Dunkl, C.: Orthogonal Polynomials of Several Variables. Cambridge
  University Press, Cambridge (2014)

\bibitem{Yger:2012}
Yger, A.: Tropical geometry and amoebas. Universit\'e Bordeaux 1, France (2000)

\bibitem{Zahabi:2021}
Zahabi, A.: Quiver asymptotics and amoeba: Instantons on toric
  $\textrm{C}$alabi-$\textrm{Y}$au divisors. Physical Review D  \textbf{103}(8)
  (2021)

\bibitem{Zhukov-Sadykov:2023}
Zhukov, T., Sadykov, T.: Computation of connected components of amoeba
  complement for polynomials of several variables. Programming and Computer
  Software  (2023), to appear

\end{thebibliography}

\end{document}